\documentclass[openright, 12pt]{article}
\usepackage{amsmath,amsfonts}
\usepackage{amsthm}
\newtheorem{theorem}{Theorem}
\newtheorem{proposition}[theorem]{Proposition}
\newtheorem{lemma}[theorem]{Lemma}
\newtheorem{remark}[theorem]{Remark}
\newtheorem{corollary}[theorem]{Corollary}
\newtheorem{example}[theorem]{Example}
\newtheorem{definition}[theorem]{Definition}
\newtheorem{conjecture}[theorem]{Conjecture}
\newtheorem{problem}[theorem]{Problem}

\usepackage{float} 
\usepackage{multirow} 
\usepackage[dvips]{graphicx}
\usepackage{psfrag}
\usepackage{enumerate}

\usepackage{tikz-cd}

\usepackage{amssymb}
\usepackage[latin1]{inputenc}

\usepackage{hyperref}

\newcommand{\iso}{{\text{Iso}}}

\newcommand{\diam}{{\text{diam}}}
\newcommand{\bet}{{\text{b}_1}}
\newcommand{\id}{{\text{Id}}}
\newcommand{\tr}{{\text{Tr}}}

\newcommand{\real}{{\text{Re}}}
\newcommand{\vol}{{\text{Vol}}}
\newcommand{\en}{{\text{End}}}
\newcommand{\kker}{{\text{Ker}}}

  \hypersetup{colorlinks}
\hypersetup{citecolor=blue}
\hypersetup{urlcolor=blue}

\begin{document}

\begin{center}

\noindent  \large Fundamental groups and group presentations with

 bounded relator lengths.\normalsize

\begin{normalsize}
Sergio Zamora

zamora@mpim-bonn.mpg.de
\end{normalsize}
\end{center}

\begin{abstract}

We study the geometry of compact geodesic spaces with trivial first Betti number admitting large finite groups of isometries. We show that if a finite group $G$ acts by isometries on a compact geodesic space $X$ whose first Betti number vanishes, then $\diam (X) / \diam (X / G ) \leq 4 \sqrt{ \vert G \vert }$.

For a group $G$ and a finite symmetric generating set $S$, $P_k(\Gamma (G, S))$ denotes the 2-dimensional CW-complex whose 1-skeleton is the Cayley graph  $\Gamma$  of $G$ with respect to $S$ and whose 2-cells are $m$-gons for $0 \leq m \leq k$, defined by the simple graph loops  of length $m$ in $\Gamma$, up to cyclic permutations.

Let $G$ be a finite abelian group with $\vert G \vert \geq 3$ and $S$ a symmetric set of generators for which $P_k(\Gamma (G,S))$ has trivial first Betti number. We show that the first nontrivial eigenvalue  $-\lambda_1$ of the Laplacian on the Cayley graph satisfies $\lambda_1 \geq 2 - 2  \cos ( 2 \pi / k ) $. We also give an explicit upper bound on the diameter of the Cayley graph of $G$ with respect to $S$ of the form $O (k^2 \vert S \vert \log \vert G \vert )$. Related explicit bounds for the Cheeger constant and Kazhdan constant of the pair $(G,S)$ are also obtained.




\end{abstract}

\section{Introduction}

\subsection{Diameter}

Compact geodesic spaces equipped with large discrete groups of isometries have been objects of great interest for a long time and several problems can be formulated in this setting \cite{benjamini-finucane-tessera,breuillard-tointon,gelander,gorodski-lange-lytchak-mendes,turing}. One natural source of such spaces are finite-sheeted Galois covers of compact Riemannian manifolds.  In 2009, Petrunin asked if one can control in an interesting way the diameter of a compact universal cover  \cite{petrunin-universal}.
\begin{problem}[Petrunin]\label{problem:anton}
\rm Let $M$ be a compact Riemannian manifold and assume it admits a compact universal cover $\tilde{M}$. What is the smallest upper bound of $\diam (\tilde{M}) / \diam (M)$ in terms of $\vert \pi_1(M) \vert $?
\end{problem}
It is not hard to show that $\diam (\tilde{M})/ \diam(M) \leq \vert \pi_1( M ) \vert $ \cite{petrunin-general}, but getting a better bound is non-trivial matter. The goal of this paper is to study this question and the global shape of compact universal covers in general.  One of our main results is the following.
\begin{theorem}\label{thm:diameter-effective} 
 \rm Let $X$ be a compact geodesic space and $G \leq \iso (X)$ a finite group of isometries. If the first Betti number $\bet (X)$ vanishes, then 
 \begin{equation*}
     \dfrac{\diam (X)}{\diam (X/ G)} \leq 4 \sqrt{ \vert G \vert }.         
 \end{equation*}    
\end{theorem}
Asymptotically as $\vert G \vert \to \infty$, there is a stronger yet non-effective bound \cite{benjamini-finucane-tessera}. 
\begin{theorem}[Benjamini--Finucane--Tessera]\label{thm:diameter-as} 
\rm Let $X_n $ be a sequence of compact geodesic spaces and $G_n \leq \iso (X_n) $ a sequence of finite groups with  $\vert G_n \vert \to \infty$ as $n\to \infty$. If the first Betti numbers $\bet (X_n)$ vanish, then for each $\varepsilon > 0$ one has
\[ \dfrac{ \diam  ( X_n  )}{ \diam  (X_n / G_n )  } = \, O \left( \vert G _n \vert ^{\varepsilon}   \right)  .  \]  
\end{theorem}
Problem \ref{problem:anton} is better handled when reformulated in terms of Cayley graphs. For a group $G$ and a finite symmetric generating set $S$ we denote by $\Gamma (G, S)$ the Cayley graph of $G$ with respect to $S$. 

For a graph $\Gamma$ and an integer $k \in \mathbb{N}$, as in \cite{delasalle-tessera} we denote by $P_k(\Gamma )$ the 2-dimensional CW-complex whose 1-skeleton is $\Gamma$ and whose 2-cells are $m$-gons for $0 \leq m \leq k$, defined by the simple graph loops  of length $m$ in $\Gamma$, up to cyclic permutations. 
\begin{proposition}[\v{S}varc--Milnor Lemma]\label{prop:sml}
  \rm  Let $X$ be a  proper geodesic space, $p \in X$,  $G \leq \iso (X)$ a discrete group, $\delta \geq 0$, and $r \geq 2 \cdot \diam (X/ G) + \delta $. Then $ S : = \{ g \in G \, \vert  \, d(gp,p) \leq r  \} $  generates $G$. Moreover, if we equip $G$ with the metric induced from $\Gamma := \Gamma (G,S)$, for all $g, h \in G$ one has
      \[  \delta \cdot  [ d_{\Gamma}(g,h) -1 ]  \leq d_X( gp , hp ) \leq r \cdot d_{\Gamma}(g,h) .   \]
\end{proposition}
\begin{proposition}\label{prop:cover}
  \rm  Let $X$ and $\Gamma$ be as in Proposition \ref{prop:sml}. Then $\pi_1(P_3(\Gamma))$ is a quotient of $\pi_1(X)$.
\end{proposition}
A proof of the \v{S}varc--Milnor Lemma can be found in \cite{delaharpe}, and Proposition \ref{prop:cover} will be proven in Section \ref{sec:coverings}.  Using these well known results, Theorem \ref{thm:diameter-effective} becomes a corollary of  its Cayley graph counterpart. 
\begin{theorem}\label{thm:diameter-effective-cayley}
 \rm  Let $k \geq 3$, $G$ be a finite group, and $S \subset G$ a finite symmetric set of generators for which $ P_k( \Gamma (G,S) )  $ has trivial first Betti number. Then 
 \begin{equation}\label{eq:diameter-effective-cayley}
   \textnormal{diam}( \Gamma (G,S) ) \leq  \left( \sqrt{4 \vert G \vert + 1}-2 \right) \left\lfloor \frac{k+2}{3} \right\rfloor   .    
 \end{equation}
 \end{theorem}
\begin{remark}\label{rem:spider}
 \rm  It is well known that for $k \geq 3$, a group $G$ and a finite symmetric set of generators $S$, the complex $P_k (\Gamma (G,S))$ is simply connected if and only if $G$ admits a presentation $\langle S \, \vert \, R \rangle $ with $R$ consisting of words of length $\leq k$  \cite[Section 2]{delasalle-tessera}\footnote{see also the primer by Yann Ollivier \href{http://www.yann-ollivier.org/maths/primer.php}{http://www.yann-ollivier.org/maths/primer.php}}. Moreover, if one considers the abstract group $\tilde{G} = \langle S \, \vert \, R_k \rangle $, where $R_k$ consists of the words of length $\leq k$ representing the identity in $G$, then $P_k (\Gamma (\tilde {G}, S))$ is the universal cover of $P_k (\Gamma (G,S))$ and the fundamental group of $P_k (\Gamma (G,S))$ is precisely the kernel of the natural map $\tilde {G} \to G$.
\end{remark}
By Remark \ref{rem:spider}, Theorem \ref{thm:diameter-effective-cayley} has the following implication.
\begin{corollary}
\rm  Let $k \geq 3$, $G$ be a finite group, and $S \subset G$ a finite symmetric set of generators for which $G$ admits a presentation $\langle S \, \vert \, R \rangle$ with $R$ consisting of words of length $\leq k$. Then (\ref{eq:diameter-effective-cayley}) holds. 
\end{corollary}

\subsection{Kazhdan constant, Cheeger constant, and spectral gap}

The \v{S}varc--Milnor Lemma implies that the medium-scale geometric features of $X$ and $\Gamma$ are closely related to each other. We now focus on such properties.  Recall that for a finite group $G$ and a finite symmetric set of generators $S$,  the Kazhdan constant $K (G,S)$, Cheeger constant $h$, and spectral gap $\lambda_1$ are related by the following inequalities
\begin{equation}\label{eq:isoperimetrics}
     \frac{h^2}{\vert S \vert ^2}  \leq   \frac{2 \lambda_1  }{\vert S \vert } \leq K(G,S) ^2 \leq 2 \lambda_1 \leq  4  h  .
 \end{equation}
We refer the reader to Section \ref{sec:isoperimetrics} for the definition of such quantities and further comments on  (\ref{eq:isoperimetrics}). For now we just mention that the three non-negative quantities $K(G,S)$, $h$, and $\lambda_1$ measure the connectivity of $\Gamma (G,S)$ in different ways. The other main result of this paper concerns finite abelian groups.
\begin{theorem}\label{thm:k-spectrum}
 \rm Let $k \geq 3$, $G$ a finite abelian group with $\vert G \vert \geq 3$, and $S\subset G$ a symmetric set of generators for which $ P_k (\Gamma (G , S))$ has trivial first Betti number. Then the Kazhdan constant satisfies
 \begin{equation}\label{eq:K}
   K(G,S)  \geq   2 \cdot \sin (  \pi /k  ).     
 \end{equation}
 Consequently, the Cheeger constant, spectral gap, and diameter satisfy
 \begin{gather}
  2  h \geq \lambda_1 \geq 2 - 2  \cos (2 \pi / k ) , \label{eq:lambda}\\
  \diam (\Gamma (G,S)) \leq  \frac{  \vert  S  \vert + 1 - \cos (2 \pi / k)   }{2 (1 -  \cos (2\pi/k)) } \log \vert G \vert + 1 . \label{eq:diameter}
  \end{gather}
  \end{theorem}  
A consequence of Theorem \ref{thm:k-spectrum} is an upper bound on the mixing time of the random walk in the corresponding Cayley graph (see Remark \ref{rem:isoperimetric-interpretation}). We refer the reader to Section \ref{sec:walks} for the definitions of random walk and mixing time. For now we just mention that $\tau_{\Gamma} (c)$ is an estimate of how long does one have to wait for heat to propagate evenly (how evenly? quantified by $c$) along the network $\Gamma$.
  \begin{corollary}\label{cor:mix}
\rm Let $k \geq 4$, $G$, and $S$ be as in Theorem \ref{thm:k-spectrum}. If 
 $\tau _{\Gamma } : [0,2] \to \mathbb{N}$ denotes the mixing time of the Cayley graph $\Gamma (G,S)$, then
\begin{equation*}
    \tau _{\Gamma} (c) \leq  \frac{k^2 \vert S \vert }{ 32} [\log \vert G \vert - 2 \log (c )] + 1 .
\end{equation*}
 \end{corollary} 
Theorem \ref{thm:k-spectrum} also yields an effective bound on the diameter of the universal cover of a closed Riemannian manifold with finite abelian fundamental group.
\begin{corollary}\label{cor:ricci}
\rm Let $M$ be a closed $n$-dimensional Riemannian manifold with $\textnormal{diam} (M)$ $ =D$, Ricci curvature $\geq \kappa (n-1)$ for some $\kappa \in \mathbb{R}$, and having a point whose injectivity radius is $\geq 2 r_0 > 0 $.  If its fundamental group  $\pi_1 (M)$ is finite and abelian, then the universal cover $\tilde{M}$ satisfies
\[  \dfrac{ \textnormal{diam}  (\tilde{M} )}{ \textnormal{diam}  (M) } \leq 4 +   \left\lfloor \left[  \frac{2 v_n^{\kappa} (2D + r_0)   }{ 3 v_n^{\kappa} (r_0)}  + \frac{1}{3}
\right] \log \vert \pi_1 (M) \vert  \right\rfloor  ,   \]
where  $v_n^{\kappa}(r)  $ denotes the volume of a ball of radius $r$ in the $n$-dimensional simply connected space of constant sectional curvature $\kappa $.
\end{corollary}
Considering the situation when $\diam (\Gamma) \to \infty$, there are bounds similar to the ones in Theorem \ref{thm:k-spectrum} for groups that are not necessarily abelian \cite{breuillard-tointon}.

\begin{theorem}[Breuillard--Tointon]\label{thm:k-spectrum-as}
\rm  Let $G_n$ be a sequence of finite groups, $S_n \subset G_n $ a sequence of finite symmetric sets of generators, and $\Gamma_n : = \Gamma (G_n, S_n)$ the corresponding Cayley graphs. Assume there is a sequence $k_n = o (\diam(\Gamma_n)) $ such that the first Betti numbers $\bet (P_{k_n}(\Gamma_n))$ vanish. Then for each $\varepsilon > 0 $, the quantities
\[    K (G_n,S_n) , \, \lambda_1 (G_n,S_n) , \,    h (G_n,S_n)  ,                             \]
cannot go to zero faster than $(\vert S _n \vert / \vert G _n \vert )^{\varepsilon}$  as $n \to \infty$.
\end{theorem}

\subsection{Outline}

In Section \ref{sec:examples} we present some computations and examples, and discuss related open problems and potential lines of research. 

In Section \ref{sec:prelim} we introduce our notation and the standard theory we will need. 

In Section \ref{sec:diameter} we give the proofs of Theorems \ref{thm:diameter-effective}, and \ref{thm:diameter-effective-cayley}. Theorem \ref{thm:diameter-effective} follows from Theorem \ref{thm:diameter-effective-cayley} which in turn depends on an elementary combinatorial argument. We also present a proof of Theorem \ref{thm:diameter-as} since it is currently stated in the literature only  in the setting of vertex-transitive graphs \cite[Theorem 1]{benjamini-finucane-tessera}. 

In Section \ref{sec:abelian} we give the proofs of Theorem \ref{thm:k-spectrum} and Corollaries \ref{cor:mix} and \ref{cor:ricci}. An elementary geometric observation yields estimate (\ref{eq:K}), from which all other results follow.

\section{Examples and further problems}\label{sec:examples}

\subsection{Diameter}

Theorem \ref{thm:diameter-as} implies that the explicit bound in Theorem \ref{thm:diameter-effective} is far from being sharp as $\vert G \vert \to \infty$. By a fundamental domain argument, even without the first Betti number assumption, one always has  
\[ \frac{\diam (X)}{ \diam (X / G)} \leq  2 \cdot  \vert G \vert   ,  \]
so Theorem \ref{thm:diameter-effective} says nothing new for $\vert G \vert \leq 4$. However, for a larger number, say, 120, Theorem \ref{thm:diameter-effective} gives a meaningful bound (again, likely far from sharp). The following example was pointed out by Kuperberg \cite{petrunin-universal}.

\begin{example}\label{example:poincare}
\rm  Let $\tilde{X} = \mathbb{S}^3$ equipped with its usual metric, and consider $X = (\tilde{X}/\sim )$ the Poincar\'e sphere \cite[Example 1.4.4 and Problem 4.4.17]{thurston}. Then $\tilde{X}$ is the 120-sheeted universal cover of $X$, and 
\[ \dfrac{\diam (\tilde{X})}{  \diam (X)} = \dfrac{\pi}{\arccos \left( \varphi^2 / \sqrt{8} \right) } \approx 8.09 ,  \]
where $\varphi $ is the golden ratio. On the other hand, the bound provided by Theorem \ref{thm:diameter-effective} is $4\sqrt{120} \approx 43.81$.
\end{example}

\begin{proof}[Proof sketch:]  A Voronoi domain of the quotient  $\tilde{X} \to X$ is a regular dodecahedron $K \subset \mathbb{S}^3$ with dihedral angles equal to $\pi/3$ \cite[Section 3.2.4]{martelli}. Let $O\in K$ be the center of the dodecahedron,  $F_1,$ $ F_2$  be the centers of two adjacent faces of $K$, and $V_1,$ $V_2$ be the vertices shared by such faces. Also let $P$ be the midpoint between $V_1$ and $V_2$. The diameter of $X$ is attained by $d_{\mathbb{S}^3}(O, V_1)$  (see Figure \ref{fig:dod}).

\begin{figure}[h]
\centering
\psfrag{a}{$O$}
\psfrag{b}{$F_2$}
\psfrag{c}{$F_1$}
\psfrag{d}{$V_2$}
\psfrag{e}{$V_1$}
\psfrag{f}{$P$}
\includegraphics[scale = 1.3]{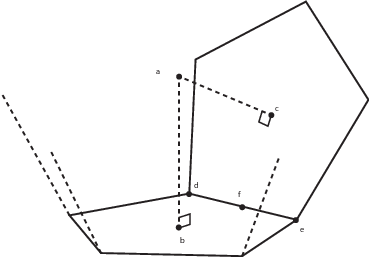}
\caption{We can use the knowledge of the angles in the triangle $OF_1P$ to deduce the length of the segment $OP$. We then proceed to compute the length of the segment $OV_1$ using the length of the segment $OP$ and the known angles of the triangle $OPV_1$.}\label{fig:dod}
\end{figure}

By the symmetry of the dodecahedron, 
\[ \measuredangle OF_1P = \measuredangle OF_2P = \measuredangle OPV_1 = \measuredangle OPV_2 =      \frac{\pi}{2}  .\]
 Using elementary geometry one can also compute the angles 
\[ \measuredangle F_1 O F_2 =  \arccos \dfrac{1}{3}  ,\text{ }  \measuredangle V_1 O V_2 =    \arccos  \dfrac{\sqrt{5}}{3}  .   \]
By the spherical laws of sines and cosines \cite{todhunter}, this is enough information to recover the length $d_{\mathbb{S}^3}(O, V_1) = \arccos ( \varphi^2 / \sqrt{8}  )$.
\end{proof}

It would also be interesting to investigate how sharp is Theorem \ref{thm:diameter-as}.  The known example in which $\diam  (X_n ) / \diam (X_n / G_n )$ grows the fastest with respect to $\vert G _n\vert $ is the following (again pointed out by Kuperberg \cite{petrunin-universal}), which naturally leads to Conjecture \ref{conj:petrunin-diameter} below.

\begin{example}
\rm  Let $G_n $ be the symmetric group (the set of bijections of the set $\{  1, \ldots , n  \}$), and $S_n$ the set of transpositions of consecutive elements of $\{1, \ldots , n    \}$ (we consider $n$ and $1$ not to be consecutive). Setting $\Gamma _ n : = \Gamma (G_n , S_n)$, we have:
\begin{enumerate}
\item $ P_6 (\Gamma_n) $ is simply connected for all $n$. \label{item:k1}
\item $\diam (\Gamma _n) = o \left(   \left(   \log \vert G _n \vert       \right)^2     \right)$.\label{item:k2}
\item  $ \left( \log \vert G _n \vert \right)^{2- \varepsilon} = o \left(  \textnormal{diam} (\Gamma_n)    \right) $ for every $\varepsilon > 0$.   \label{item:k3}
\end{enumerate}
\end{example}

\begin{proof}[Proof sketch:] \ref{item:k1}: The group $G_n $ can be presented as $\langle  S_n \vert R_n   \rangle$, with $S_n=\{\sigma _1 , \ldots  ,$ $ \sigma_{n-1}    \}$, and  $R_n $ consisting of the words $\sigma_i^2$ for all $i$, $(\sigma_i\sigma_{i+1})^3$ for all $i=1, \ldots , n-2$, and $(\sigma_i  \sigma_j)^2$ with $\vert i-j \vert \geq 2$.  Since each word in $R_n$ has length $\leq 6$, the complex $P_6(\Gamma_n)$ is simply connected. 

\ref{item:k2} and \ref{item:k3}: Every permutation in $G_n$ can be written as a composition of at most $n(n-1)/2$ elements in $S_n$, where the maximum is achieved by the permutation 
\[ i\to( n+1-i) , \, i \in \{1, \ldots , n \} \]
 that ``reverses'' the order. This means that 
 \[ \diam (\Gamma_n) = n (n-1)/2,\]
  while 
  \[ \log \vert G_n \vert = \log (n !) ,\]
   which is of the order of $n \log n$. 
\end{proof}

\begin{conjecture}[Petrunin]\label{conj:petrunin-diameter}
\rm There is  $C>0$ such that if $G$ is a finite group and $S$ a set of generators for which $P_3(\Gamma (G,S))$ is simply connected,  then 
\[          \diam (\Gamma (G,S)) = O \left( (\log \vert G \vert )^C \right).                  \]
\end{conjecture}
This question draws resemblance to another well known problem \cite{babai-seress}.
\begin{conjecture}[Babai]\label{conj:babai}
 \rm There is $C>0$ such that if $G$ is a finite non-abelian simple group and $S\subset G$ is any set of generators, then
 \[      \diam (  \Gamma (G, S) ) =   O \left( (\log \vert G \vert )^C \right).                   \]
\end{conjecture}
Note however, that Babai's Conjecture concerns any set of generators, while Petrunin's Conjecture is about geometrically chosen sets of generators. It would be interesting to investigate how intertwined these two problems are. For instance, does the hypothesis in Conjecture \ref{conj:petrunin-diameter} of $P_3(\Gamma (G,S))$ being simply connected imply that the graph $\Gamma (G,S) $ looks like the Cayley graph of a non-abelian simple group? We refer the reader to  \cite{eberhard-jezernik,halasi} for recent updates on the state of Conjecture \ref{conj:babai}. 

We would like to also point out that Conjecture \ref{conj:petrunin-diameter} is still very interesting when the group $G$ is abelian, in which case $C$ could even be 1. For abelian groups, the known example in which $\diam (\Gamma (G,S ))$ grows the fastest with respect to $\vert G \vert $ is the following (pointed out by Petrunin \cite{petrunin-universal}).
\begin{example}
\rm  Let $G_n = \mathbb{Z} /  ( 2^n) \mathbb{Z} $, and 
\[ S_n = \{    \pm  1 , \pm 2, \pm 2^2 ,  \ldots , \pm 2^{n-1} \}.\]
Setting $\Gamma_n : = \Gamma (G_n, S_n)$, we have:
\begin{enumerate}
\item $P_3(\Gamma_n)$ is simply connected for all $n$.\label{item:p1}
\item $\diam (\Gamma _n) = O \left(     \log \vert G _n \vert     \right)$.\label{item:p2}
\item  $\log \vert G _n \vert = O \left(  \textnormal{diam} (\Gamma_n)    \right).  $   \label{item:p3}
\end{enumerate}
\end{example}

\begin{proof}[Proof sketch:] \ref{item:p1}: The group $G_n$ can be presented as $\langle S_n \vert R_n \rangle$, where $R_n$ consists of the expressions $ 2^j - 2^{j-1} - 2^{j-1} =0$ for $j\in \{ 1, \ldots , n-1 \}$, and $2^{n-1} + 2^{n-1} = 0$. Since each word in $R_n$ has length $\leq 3$, the complex $P_3(\Gamma_n)$ is simply connected. 

\ref{item:p2}: Any number in $\{ 1, 2 , 3, \ldots , 2^n -1 \}$ can be written (using binary base) as a sum of at most $n$ summands of the form $2^j$, $j \in \{1, \ldots , n-1\}$. Hence 
\[ \text{diam}(\Gamma_n) = O (n) = O(\log \vert G_n \vert).\]

\ref{item:p3}: Given a sequence of length $n$ of $0$'s and $1$'s, one could count the number of ``jumps'' from one digit to another. E.g., 0001111 has 1 jump, 0011101 has 3 jumps, 1010110 has 5 jumps, etc. By writing an element $x \in G_n$ in binary base, we obtain a sequence of length $n$ of $0$'s and $1$'s. One can check that the effect of adding or substracting a power of $2$ to $x$ increases the number of such jumps by at most 2. 

Expressing $x= 1 + 2^2 + \ldots + 2^{2 \lfloor \frac{n-1}{2} \rfloor }$ in binary we find $n-1$ jumps. Hence at least $\left\lfloor \frac{n}{2} \right\rfloor$ elements of $S_n$ are required to write down $x$. This implies,
\[  \log \vert G_n \vert = O \left (  \left \lfloor   \dfrac{n}{2}    \right\rfloor  \right) = O (\text{diam}(\Gamma_n))  .  \]
\end{proof}

We conclude this section by pointing out that examples of spaces $X$ with $\bet (X) = 0$ and finite groups $G \leq \iso (X)$ with $\diam (X) / \diam (X / G )$ of order $\log \vert G \vert $ arise in number theory.

\begin{example}[Calegari--Dunfield, Boston--Ellenberg]\label{exam:hyperbolic}
\rm  There is a closed hyperbolic 3-manifold $M_0$ admitting a tower of regular $m_n$-sheeted covers $M_n \to M_0$ with $\bet (M_n) = 0$ for all $n$, and  
\begin{equation}\label{eq:brooks-diameter}
    c \cdot \log (m_n) \leq   \frac{ \diam (M_n ) }{\diam (M_0)}  \leq C \cdot \log (m_n)  
\end{equation}   
for some $C > c > 0$.
\end{example}

\begin{proof}[Proof sketch:] In \cite[Theorem 1.6]{calegari-dunfield}, a sequence of regular finite-sheeted covers $M_n \to M_0$ is  constructed, with $M_0$ a quotient of the hyperbolic space $\mathbb{H}^3$ by an arithmetic lattice $\Gamma  \leq PGL_2(\mathbb{C})$, which possesses the Selberg property  (see \cite[Section 2]{lubotzky-i}). 

Since these covers correspond to congruence subgroups of $\Gamma$, by the Selberg property their (analytic) spectral gaps satisfy $\lambda_1 (M_n) \geq \varepsilon $ for some  $\varepsilon > 0$. Then by the work of Brooks \cite[Theorem 1]{brooks-iii}, the estimates (\ref{eq:brooks-diameter}) follow.

In \cite{boston-ellenberg} it was then  proven that $\bet (M_n) = 0$ for all $n$ (this fact was proven initially in \cite{calegari-dunfield} assuming the Generalized Riemann Hypothesis and Langlands-type conjectures).
\end{proof}

\subsection{Spectral gap}

The bounds in Theorem \ref{thm:k-spectrum} and Corollary \ref{cor:mix} are rather general, so we don't expect them to be fully sharp.

\begin{example}
\rm Let $G_n = (\mathbb{Z}/ 2 \mathbb{Z})^n$,  $S_n  = \{ e_i \in G_n \vert i \in \{ 1, \ldots , n \} \} $, where $e_i$ is the $i$-th basis vector, and $\Gamma_n  = \Gamma (G_n, S_n)$. Then $P_4(\Gamma_n)$ is simply connected for each $n$ so by Corollary \ref{cor:mix}
\[  \tau_{\Gamma_n } (c) \leq   \dfrac{n }{2}   [ n  \log 2 - 2 \log (c)  ] + 1 .     \]
For fixed $c>0$, this bound grows quadratically in $n$. However, a careful computation  \cite{diaconis} shows that
\[    \tau_{\Gamma_n}(c)  = O ( n \log n  ) .  \]
\end{example}
\begin{example}
\rm Let $G_n = \mathbb{Z}/n\mathbb{Z}$,  $S_n = \{   -1,1 \}  $, and $\Gamma_n = \Gamma (G_n, S_n)$.  For fixed $\varepsilon >0$, consider the sequence  
\[  t_n := \left\lfloor n^{2-\varepsilon }\right\rfloor .    \]
 By the Central Limit Theorem \cite[Section 7.3]{ash}, the random walk in $\Gamma_n$ after $t_n$ steps will be concentrated in the interval $[-n/\sqrt{t_n}, n/\sqrt{t_n}]$. That is, 
\[ W^{t_n}\left(   \left[ - \frac{n}{\sqrt{t_n}},\frac{ n}{\sqrt{t_n}} \right]  \right) \to 1 \text{ as }n \to \infty .   \]
 This implies that for fixed $c <2 $, and large enough $n$,
\[  \tau_{\Gamma_n}(c)  \geq  t_n \geq  n^{2 - \varepsilon}-1  . \]
Notice $P_n(\Gamma_n)$ is simply connected for each $n$, so the bound given by Corollary \ref{cor:mix} is the following,  not far from being sharp:
\[  \tau_{\Gamma_n} (c) \leq     \dfrac{ n^2 }{16}   [ \log n - 2 \log (c)  ] + 1 .            \]
\end{example}

It would be desirable to find explicit bounds similar to the ones of Theorem \ref{thm:k-spectrum} in the non-abelian setting. However, at the moment the topological condition of $P_k(\Gamma(G,S))$ having trivial first Betti number (or even being simply connected) seems very hard to use when studying isometric actions $G \to \iso (\mathbb{S}^n)$ with $n \geq 2$. In this direction, there is a universal control on the diameters of quotients of spheres by group actions \cite{gorodski-lange-lytchak-mendes}. 

\begin{theorem}[Gorodski--Lange--Lytchak--Mendes]
    \rm There is $\delta > 0 $ such that for any $n \geq 2$ and any compact group $G \leq \iso (\mathbb{S}^n)$ not acting transitively, one has
    \[   \diam (\mathbb{S}^n / G ) \geq \delta .       \]
\end{theorem}
With techniques similar to the ones in the work of Mantuano \cite{mantuano}, it seems possible to recover, using Theorem \ref{thm:k-spectrum}, effective estimates on spectral gaps and medium-scale isoperimetric inequalities for compact Riemannian manifolds with trivial first Betti number and actions by finite abelian groups with small quotient. Successful results in similar programs have been obtained by Brooks \cite{brooks-i,brooks-ii}, Buser  \cite{buser}, Burger \cite{burger},  Magee \cite{magee},  and several others, mostly for surfaces.

\subsection{(Lack of) Gromov--Hausdorff precompactness}

An interesting problem in the theory of finite groups was to understand the possible limits of finite homogeneous spaces. For instance; can one find a sequence of compact geodesic spaces $X_n$ and finite groups $G_n \leq \iso (X_n)$ with $\diam (X_n / G_n )\to 0$ such that $X_n$ converges to $\mathbb{S}^2$ in the Gromov--Hausdorff sense? This question was answered negatively by Turing \cite{turing}, and building upon his work, Gelander \cite{gelander} proved the following.
\begin{theorem}[Gelander]\label{thm:gelander}
 \rm Let $X_n$ be a sequence of compact geodesic spaces and $G_n \leq \iso (X_n)$ a sequence of finite groups with $\diam (X_n  /G_n ) \to 0$. If $X_n$ converges in the Gromov--Hausdorff sense to a compact space $X$, then $X$ is a (possibly infinite-dimensional) torus.
\end{theorem}
A consequence of Theorem \ref{thm:gelander} is that a sequence of normalized universal covers cannot have a ``limit shape''.
\begin{corollary}\label{cor:gh}
\rm Let $X_n$ be a sequence of compact geodesic spaces and $G_n \leq \iso (X_n)$ a sequence of finite groups with $\diam(X_n/G_n) / \diam (X_n) \to  0 $. If  $\bet (X_n) = 0$ for all $n$, then the sequence $X_n / \diam (X_n)$ diverges in the Gromov--Hausdorff sense.  
\end{corollary}
\begin{proof}
 Assuming the contrary,  $X_n / \diam (X_n)$ converges to a space $X$ of diameter $1$. By Theorem \ref{thm:gelander}, $X$ is a torus so it admits a regular covering with Galois group $\mathbb{Z}$. Then by   the work of Sormani--Wei \cite[Theorem 3.4]{sormani-wei},  there are surjective morphisms $\pi_1 (X_n) \to \mathbb{Z}$ for $n$ large enough contradicting the assumption $\bet (X_n) = 0$.  
\end{proof}

\begin{remark}
\rm Theorems \ref{thm:diameter-as} and \ref{thm:k-spectrum-as} are proven in a similar fashion. In \cite{benjamini-finucane-tessera,breuillard-tointon}, building upon the structure of approximate groups by Breuillard--Green--Tao \cite{breuillard-green-tao}, it is proven that if one had contradicting subsequences, then the normalized spaces would converge to a finite-dimensional torus, contradicting the lower-semi-continuity of the first Betti number \cite{sormani-wei}. 
\end{remark}

It would be interesting to further understand what causes the behavior of the sequences $X_n / \diam (X_n)$ in Corollary \ref{cor:gh}. Recall that some known families of Gromov--Hausdorff divergent sequences such as $X_n = \mathbb{S}^n$ or $X_n = (\mathbb{Z}/2\mathbb{Z})^n$ present a concentration of measure property \cite{ledoux}. 

\begin{definition}
\rm We say that a sequence $(X_n, d_n , \mu _n )$ of metric probability spaces of diameter $1$ is a \textit{Levy family} if for any sequence of $1$-Lipschitz maps $f_n : X_n \to \mathbb{R}$, the sequence of variances $    \text{Var} \left(  (f_n )_{\ast} \mu_n \right)    $ goes to $0$ as $n \to \infty$.
\end{definition}
 
\begin{conjecture}[Petrunin]
\rm Let $(X_n, d_n, \mu_n )$ be a sequence of compact simply connected geodesic probability spaces and $G_n \leq \iso (X_n)$ a sequence of finite groups of measure preserving isometries with $\diam (X_n / G_n) \to 0$. Then $X_n$ is a Levy family.  
\end{conjecture}

\section{Preliminaries}\label{sec:prelim}

\subsection{Notation}

For a finite-dimensional $\mathbb{C}$-Hilbert space $V$, we denote by $ \en (V)$ the space of linear maps $V \to V$ and by $U (V) \subset \en (V) $ the set of unitary automorphisms. If $V = \mathbb{C}^n$, then we denote $U (V)$ also by $U (n)$. For $A \in \en (V) $, we denote its spectrum by $\sigma (A) \subset \mathbb{C}$ and its adjoint by $A^{\ast}$. The trace operator is denoted by $ \tr: \en (V ) \to \mathbb{C}  $. When $V$ is $1$-dimensional, we will identify $\en (V )$   with $\mathbb{C}$  via $\tr: \en (V ) \to \mathbb{C}$.

For a path connected topological space $X$, we denote its first Betti number by $\bet (X)$. Recall that it equals the supremum of the $m$ for which there is a surjective morphism $\pi_1(X) \to \mathbb{Z}^m$.  

For a metric space $X$, $p \in X$, and $r > 0 $, we denote by $B(p,r)$ the open ball of radius $r$ around $p$. For two metric spaces $X$ and $Y$, we denote their Gromov--Hausdorff distance by $d_{GH}(X,Y)$.

\subsection{Graphs and CW-complexes}\label{sec:graphs}

For the purposes of this paper, a graph always means a locally finite undirected graph without loops or multiple edges. For vertices $x,y$ in a graph, we write $x \sim y$ if there is an edge connecting $x$ to $y$. For an edge $[x,y]$, we denote by $(x,y)$ its interior. 

For a sequence of vertices $v_0, v_1, \ldots , v_m$ in a graph $\Gamma$ such that $v_{i-1} \sim v_i$ for each $i \in \{ 1, \ldots , m \}$,  we denote by $[ v_0 , \ldots , v_m ] $ the curve $\gamma: [0,m] \to \Gamma$ with $\gamma (i) = v_i$ for every $i \in \{ 0, \ldots , m \}$, so that $\gamma  | _{[i-1,i]} $ travels along the edge $v_{i-1}v_i$.  A curve (loop) of this form is called a \textit{graph curve (loop)} of length $m$.

For vertices $x,y $ in a connected graph $\Gamma$, the \textit{graph distance} $d_{\Gamma}(x,y)$ between $x$ and $y$ is the minimum $m$ for which there is a graph curve of length $m$ connecting them.

For a graph $\Gamma$ and an integer $k \in \mathbb{N}$,  we denote by $P_k(\Gamma )$ the 2-dimensional CW-complex whose 1-skeleton is $\Gamma$ and whose 2-cells are $m$-gons for $0 \leq m \leq k$, defined by the simple graph loops of length $m$ in $\Gamma$, up to cyclic permutations. 
\begin{remark}\label{rem:p-metric}
\rm It is not hard to equip $P_k(\Gamma)$ with a geodesic metric that restricted to $\Gamma $ coincides with its original metric, and such that $d_{GH}(P_k(\Gamma), \Gamma ) \leq k $. For instance; for $k=3$ one can make each $2$-cell a Reuleaux triangle. 
\end{remark}
Let $G $ be a group and $S \subset G$ a symmetric generating subset. The \textit{Cayley graph} $\Gamma (G,S)$ of $G$ with respect to $S$ is defined to be the one with $G$ as its vertex set and such that two distinct elements $g, h \in G$ are adjacent  if and only if $g = hs$ for some $s \in S$. 

We now state a trivial observation. We include its proof since this same counting argument will be used later (see Claim 2 in the proof of Theorem \ref{thm:diameter-effective-cayley}).

\begin{lemma}\label{lem:t-not-separating}
  \rm   Let $G$ be a finite group with $\vert G \vert \geq 3$, $S \subset G$ a symmetric set of generators, $\Gamma : = \Gamma (G, S)$ the Cayley graph, and  $t \in S \backslash \{ e \} $. Then $\Gamma \backslash  (e,t ) $ is connected.
\end{lemma}

\begin{proof}
    If $t$ has order $m \geq 3$, then the path $[t,t^2, \ldots , t^m ]$ connects the endpoints of the removed edge, so we can assume $t =t^{-1}$. Let $C_1 $ and $C_2$ denote the connected components of $\Gamma \backslash (e,t)$ containing $e$ and $t$, respectively. Since multiplication by $t$ exchanges $e$ and $t$, it sends $C_1$ to $C_2$ and vice-versa, so $\vert C_1 \vert  = \vert C_2 \vert$. 

    Since $\vert G \vert \geq 3$, $S \backslash \{ e \}$ contains an element $s \neq t$. Multiplication by $s$ sends $(e,t)$ to $(s,st)$, so $s C_1$ is the connected component of $\Gamma \backslash (s,st) $ containing $s$. Since $s \neq t = t^{-1}$, the three segments $(s,e)$, $ (s,st)$, and $ (e,t)$ are distinct. Hence, 
    \begin{itemize}
        \item the path $[s,e,t]$ lies entirely in $sC_1\subset \Gamma \backslash (s,st)$.
        \item the path $[e,s,st]$ lies entirely in $C_1 \subset \Gamma \backslash (e,t)$. 
    \end{itemize}
 If $C_1 \neq C_2$, then $C_1 \cap C_2 = \emptyset$ and the above implies that 
 \begin{itemize}
     \item the connected set $ [s,e,t]  \cup C_2$ lies entirely in $sC_1 \subset  \Gamma \backslash (s,st)$.
     \item  $\{ s, e \} \cap C_2 = \emptyset$.
 \end{itemize}
 Therefore
 \[  \vert C_1 \vert = \vert sC_1 \vert \geq \vert [s, e,t] \cup C_2  \vert = 2 +  \vert C_2 \vert  =  2 + \vert C_1 \vert .               \]
    This contradiction finishes the proof of Lemma \ref{lem:t-not-separating}.
\end{proof}

\subsection{Constructing covering spaces}\label{sec:coverings}

In this section we prove Proposition \ref{prop:cover}. In order to do so, we present a general construction (cf. \cite[Section 5D]{gromov}).

\begin{proposition}\label{prop:monodromy}
    \rm  Let $X$ be a proper geodesic space, $p \in X$,   $G \leq \iso (X)$   a discrete group of isometries, and $r \geq 2 \cdot \diam (X/ G)$. Then set $S: = \{ g \in G \vert d(gp,p ) \leq r \}$, and let $\tilde{G}$ be the abstract group  generated by $S$, with relations 
\begin{center}
$s= s_1 s_2$ in $ \tilde{G} $, whenever $s, s_1, s_2 \in S$ and $s= s_1s_2$ in $G$. 
\end{center}
Denote the  canonical embedding $S  \hookrightarrow \tilde{G} $ as $(s \to s^{\sharp } )$, and by $\Phi : \tilde{G} \to G$ the unique morphism with $\Phi (s^{\sharp }) = s$ for all $s \in S$. Then there is a regular covering $\tilde{X} \to X$ with Galois group $\kker (\Phi)$.
\end{proposition}
\begin{proof}
In order to construct the space $\tilde{X}$, notice that by discreteness of $G$, there is $\eta > 0 $ with $ S = \{ g \in G \vert d(gp,p ) < r + 2 \eta \} $. Set $B  : = B(p, r/2 + \eta ) $. Then $S = \{ g \in G \vert B \cap gB \neq \emptyset \}$. Equip $ \tilde{G} $ with the discrete topology, and consider the topological space
\begin{center}
$\tilde{X} := \left(  \tilde{G }  \times B \right) / \sim$,
\end{center}
where $\sim$ is the minimal equivalence relation such that 
\begin{equation}\label{eq:x-tilde-def}
 (gs^{\sharp}, x) \sim (g,s x)   \text{ whenever }s \in S, \, x , sx \in B. 
\end{equation}
We then obtain a continuous map $\Psi : \tilde{X} \to X$ given by
\begin{center}
$\Psi  (g, x) : = \Phi (g)(x)      .$
\end{center}
Fix $g _0 \in \tilde{G} $ and set $U : = \Phi (g_0) (B)$.  The proof of \cite[Theorem 2.32]{zamora} carries over (with $V = B$ and $\Gamma = G$) to show that $U$ is evenly covered. As $g_0$ ranges over $\tilde{G}$, the sets $\Phi (g_0)(B)$ cover $X$, so $\Psi $ is a covering map. The proof of \cite[Theorem 2.32]{zamora} again carries over to show that $\Psi$ is regular with Galois group $\kker (\Phi )$.
\end{proof} 
 
We now prove Proposition \ref{prop:cover}. Let $X, G, S, \Gamma $ be as in the statement of the proposition. Let $\tilde{G}$ be the group with presentation $ \langle S \,  \vert  \, R \rangle$, where $R$ consists of the words of length $\leq 3$ that represent the trivial element of $G$.   Then $P_3(\Gamma (\tilde{G}, S))$ is the universal cover of $P_3(\Gamma )$, and $\pi_1(P_3(\Gamma ))$ is isomorphic to the kernel of the natural map $\tilde{G} \to G$  (see Remark \ref{rem:spider}). By Proposition \ref{prop:monodromy}, there is a regular covering map $\tilde{X} \to X$ with Galois group $\pi_1(P_3(\Gamma ))$, so there is a surjective map $\pi_1(X) \to \pi_1(P_3(\Gamma ))$.

\subsection{Representation theory of finite groups}
In this section we recall the results from representation theory we will need. We refer the reader to \cite[Chapters 1-2]{serre} for proofs and further discussion.  Throughout this section, let $G$ be a finite group. 

For our purposes, a \textit{(unitary) representation} is a morphism $\rho : G \to U (V) $ for some finite-dimensional $\mathbb{C}$-Hilbert space $V$. The dimension of $V$ is called the \textit{dimension} of the representation and will be denoted by $d_{\rho}$.  We say that such representation is \textit{irreducible} if whenever there is a subspace $W \leq V$ invariant under the $G$-action, either $W = \{ 0 \}$ or $W = V$. The representation $\rho$ is said to be \textit{trivial} if $\rho (g) = \id _V$ for all $g \in G$.

Given two representations $\rho _ 1 : G \to U (V_1)$, $\rho_2 : G \to U (V_2)$, we say a linear map  $\lambda : V_1 \to V_2$  is \textit{equivariant} if 
\[      \lambda \rho_1 (g) = \rho _2 (g) \lambda \, \text{ for all } \, g \in G       .  \]
We say that $\rho_1$ and $\rho_2$ are  \textit{isomorphic} if there is an equivariant linear isomorphism $\lambda : V_1 \to V_2$. It turns out there are only finitely many isomorphism classes of irreducible representations and they satisfy
\begin{equation}\label{eq:sum-dimensions}
    \vert G \vert   = \sum _{\rho } d_{\rho}^2,
\end{equation}
where the sum is taken among the isomorphism classes of irreducible unitary representations of $G$.

Let $\mathbb{C}[G]$ be the space of functions $G \to \mathbb{C}$.  For $f \in \mathbb{C}[G]$, we denote by $f^{\ast} \in \mathbb{C}[G] $ the function given by $f ^{\ast}(g) : = \overline{f(g^{-1})}$. For $f, h \in \mathbb{C}[G]$, their product is defined as
\[   ( f \ast h )(g) : =   \sum _{u \in G} f(gu)h(u^{-1}).                \]
The  Hermitian product 
\[     \langle f,  h \rangle : = \frac{1}{\vert G \vert } \sum_{g \in G}  f(g) \overline{h(g)}                                \]
makes $\mathbb{C}[G]$ a Hilbert space called the \textit{convolution algebra} of $G$. It admits a natural action $\rho_0 : G \to U ( \mathbb{C}[G])$ given by 
\[ \rho_0(g)(f)(h) : = f ( g^{-1}h)\, \text{ for }\, f \in \mathbb{C}[G],\,  g,h \in G . \]
For $f \in \mathbb{C}[G]$ and a representation $\rho : G \to U(V)$, the \textit{Fourier transform} $\hat{f}(\rho) \in \en (V)$ is defined as
\[     \hat{f} (\rho ) : = \sum_{g \in G} f(g) \rho (g).       \]
The Fourier transform is compatible with products and adjoints. That is, 
\begin{equation}\label{eq:fourier-properties}
 \widehat{f\ast h} (\rho) = \hat{f}(\rho) \hat{h}(\rho)     \,\,\,\,\,\,\,\,\,\,\,\, \,\,\, \,\,\,\,\,  \widehat{f^{\ast} } (\rho) = \left[ \hat{f}(\rho) \right] ^{\ast}       \end{equation}
for any $f, h \in \mathbb{C}[G]$ and any representation $\rho$.

It can be shown that for any irreducible representation $\rho : G \to U (V)$, there is a (non-unique!) equivariant isometric embedding $ \iota_V: V \to \mathbb{C} [G]$. Moreover, such embeddings among all irreducible representations span $\mathbb{C}[G]$. This leads to the following important result.
\begin{theorem}[Plancherel formula]
   \rm  For $f, h \in \mathbb{C}[G]$ one has
    \[       \langle f , h^{\ast} \rangle = \frac{1}{\vert G \vert ^2 } \sum_{\rho} d_{\rho} \, \tr  ( \hat{f}(\rho) \hat{h} (\rho)   ),                                    \]
    where the sum is taken among the isomorphism classes of irreducible unitary representations of $G$.
\end{theorem}
The following result follows from the fact that any unitary representation is totally reducible.
\begin{lemma}\label{lem:x-irreducible}
 \rm For any unitary representation $\rho : G \to U (V)$, there is $x \in V$ with $\vert x \vert = 1 $ such that the span of the $G$-orbit of $x$ is irreducible. 
\end{lemma}
The ensuing result follows from the fact that commuting unitary automorphisms are simultaneously diagonalizable. 
\begin{proposition}\label{prop:abelian-1}
    \rm If $G$ is abelian, then any irreducible unitary representation is $1$-dimensional.
\end{proposition}

\subsection{Kazhdan constant, Cheeger constant, and spectral gap}\label{sec:isoperimetrics}

We refer the reader to \cite{bekka-delaharpe-valette} for a detailed introduction to the theory of Kazhdan's property (T) and related topics. Throughout this section, let $G$ be a finite group, $S \subset G$ a symmetric generating set, and $\Gamma : = \Gamma (G,S)$ the corresponding Cayley graph. 

 The \textit{Kazhdan constant} of $G$ with respect to $S$ is defined as
\[    K(G,S) : = \inf_{\rho} \inf_{\vert x \vert = 1 } \sup _{s \in S} d(\rho (s)x,x),           \]
where the first infimum is taken among non-trivial irreducible unitary representations  $\rho : G \to U (V)$ and the second infimum  is taken among unit vectors $x \in V$. 

The \textit{Laplacian} is the map $\Delta : \mathbb{C}[G] \to \mathbb{C}[G]  $ defined by
\[  \Delta (f) (g) : = \sum_{s \in S} [ f(gs) - f(g) ].                   \]
This map is equivariant, self adjoint, and its spectrum is a finite set of non-positive real numbers 
\[       \sigma (\Delta ) : = \{  0 = \lambda_0 >  - \lambda_1 \geq \ldots \geq  - \lambda_{\vert G \vert -1 }  \} .  \]
We denote the $\Delta$-eigenspace of an eigenvalue $- \lambda$ by $E_{\lambda}$. The quantity $\lambda_1$ is called the \textit{spectral gap} of $G$ with respect to $S$.

\begin{figure}[h!]
\centering
\psfrag{a}{$x$}
\psfrag{b}{$\rho (s) x$}
\psfrag{c}{$2  \sin (\theta /2 )$}
\psfrag{d}{$\theta$}
\psfrag{e}{$\cos (\theta )$}
\includegraphics[scale = 1.1]{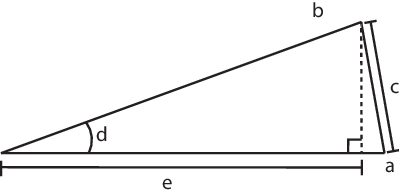}
\caption{Basic trigonometry shows that if the angle between the unit vectors $x$ and $\rho (s) x$  is $\theta$, then the distance between the endpoints is $2 \sin (\theta / 2)$.}\label{fig:trig}
\end{figure}

The \textit{Cheeger (isoperimetric) constant} of $G$ with respect to $S$ is defined by
\[  h(G,S) : =  \inf \left\{  \frac{\vert \partial A \vert }{ \vert A \vert } \, \Huge{\vert} \, A \subset G, \,  \vert A \vert \leq \frac{1}{2} \vert G \vert \right\}      ,       \]
where $\partial A$ denotes the set of edges in $\Gamma  $ connecting a vertex in $A$ with a vertex in $G \backslash A$. 
\begin{remark}\label{rem:isoperimetric-interpretation}
  \rm  The Kazhdan constant $K (G,S)$ quantifies how easy it is to tell apart isometric actions on spheres with fixed points from isometric actions on spheres without fixed points, the spectral gap $\lambda_1$ quantifies how fast heat flows through $\Gamma$, and $h$ quantifies how bad the bottlenecks of $\Gamma $ are. Each quantity measures in some way how robust is the network $\Gamma $.
\end{remark}
Recall that these quantities satisfy the well known relations
\begin{equation}\label{eq:isoperimetrics-2}
     \frac{h^2}{\vert S \vert ^2}  \leq   \frac{2 \lambda_1  }{\vert S \vert } \leq K(G,S) ^2 \leq 2 \lambda_1 \leq  4   h  
 \end{equation}
(see \cite[Section 1.2]{breuillard-tointon} for a similar expression involving $\diam (\Gamma )$).  A proof of the first and last inequalities, known as the discrete Cheeger--Buser inequalities, along with historical background can be found in \cite[Section 4.2]{lubotzky}. To verify the second inequality, take an arbitrary irreducible unitary representation $\rho : G \to U (V)$. Since there exists an equivariant embedding $\iota _V : V \to \mathbb{C}[G]$, we can assume $V \leq \mathbb{C} [V]$. Moreover, since $\Delta$ is equivariant, we can further assume $V \leq E_{\lambda}$ for some $\lambda \neq 0$. Then assume a unit vector $x \in V$ satisfies
\[    \sup_{s \in S} d(\rho_0(s) x ,x ) \leq 2 \cdot \sin ( \theta / 2 )   \]
for some $\theta \in [ 0 , \pi ] $. Then the angle between $x$ and $\rho _0 (s) x$ is at most $ \theta $ for each $s \in S$ (see Figure \ref{fig:trig}). This implies 
\[ 
  \lambda =    \langle -  \Delta x , x \rangle  = \sum _{s \in S} \langle x -  \rho _0 (s) x  ,x \rangle  \leq \vert S \vert [ 1- \cos ( \theta )  ].
\]
From the identity $ 1 - \cos (\theta)  =  2 \cdot \sin^2 (\theta / 2) $, we deduce $[ 2 \cdot \sin (\theta / 2 )]^2 \geq 2 \lambda / \vert S \vert $. Since $x$ was arbitrary and $\lambda \geq \lambda_1$, the second inequality of (\ref{eq:isoperimetrics-2}) follows. 

To verify the third inequality of (\ref{eq:isoperimetrics-2}), recall that by Lemma \ref{lem:x-irreducible}, there is a unit vector $x \in E_{\lambda_1}$ for which its $G$-orbit spans an irreducible representation. By  definition there is $t \in S$ with $d(\rho_0(t)x,x) \geq K (G,S)$. Set $K (G,S) = 2 \cdot \sin (\theta /2)$ with $\theta \in [0, \pi ]$. Then since $\vert ( \langle \rho (s) x, x\rangle \vert \leq 1 $ for all $s \in S \backslash \{ t \}$, one has 
\[    \lambda _1 = \langle - \Delta x, x \rangle = \sum_{s \in S} \langle x - \rho _0(s) x  , x \rangle \geq  1 -  \cos (\theta) =  \frac{K (G,S)^2}{2} .          \]


\subsection{Random walks on Cayley graphs}\label{sec:walks}

We refer to  \cite[Chapter 3]{diaconis} for an introduction on the theory of random walks in finite groups. Throughout this section, let $G$ be a finite group, $S \subset G$ a symmetric generating set, and $\Gamma : = \Gamma (G,S)$ the corresponding Cayley graph.

For $\alpha \geq \vert S \vert $, the \textit{random walk} on $\Gamma $ is the $G$-valued Markov process $\{   W^t_{\alpha}\}_{t \in \mathbb{N}} $ such that $W^0_{\alpha} \equiv e$, and at each time, if the walker is at $g \in G$, then it stays at $g$ with probability $1 - \vert S \vert / \alpha $ and jumps to a neighbor uniformly at random with probability $\vert S \vert / \alpha $.  This gives rise to the law
\[ \mathbb{P}\left[ W^{t+1}_{\alpha}  = g    \right] =  \frac{1}{\alpha} \left[ ( \alpha - \vert S \vert ) \mathbb{P} \left[ W^{t}_{\alpha}  = g     \right]    + \sum_{s \in S}  \mathbb{P} \left[ W^{t} _{\alpha} =  gs    \right] \right] .  \]
We will denote $W_{\alpha}^1$ simply by $W_{\alpha}$. We can identify $\mathbb{C}[G]$ with the set of complex valued measures on $G$ via the correspondence
\[  \mu \in \mathbb{C}[G] \longleftrightarrow  \mu (A) := \sum_{g \in A} \mu (g).       \]
Note that after this identification, 
\begin{equation}\label{eq:w1-expression}
    W_{\alpha} = \delta_{e} + \frac{1}{\alpha } \Delta (\delta _e ) ,   
\end{equation}
where $\delta_e$ denotes the Dirac mass at $e$. It is also straightforward to verify that 
\begin{equation} \label{eq:w-properties}
W^s_{\alpha} \ast W^t_{\alpha} = W^{s+t}_{\alpha}  \,\,\,\,\,\,\,\,\,\,\,\,\, \left(  W^t_{\alpha} \right)^{\ast} = W^t_{\alpha}   
\end{equation}
 for all $s,t \in \mathbb{N}$. If $\alpha > \vert S \vert $, the distribution $W^t_{\alpha}$ converges to the uniform distribution $U $ on $G$ as $t \to \infty$. One can quantify how fast this convergence occurs with the quantity
\[  \varepsilon _{\alpha }(t)  : =  \sum_{g \in G}  \Huge\lvert W^t_{\alpha} (g) - \frac{1}{\vert G \vert } \Huge\rvert  .    \] 
The \textit{mixing time} $\tau _{\Gamma }^{\alpha}: [0,2 ] \to \mathbb{N}$ of the process $\{  W^t_{\alpha}   \}_{t\in \mathbb{N}}$  is defined as 
\[  \tau_{\Gamma }^{\alpha} (c) : = \inf \{ t \in \mathbb{N} \, \vert \, \varepsilon _{\alpha} (t) \leq c     \}     .\]
A direct computation shows that for all $t \in \mathbb{N}$,  
\begin{equation}\label{eq:w-u-product}
     \langle W_{\alpha} ^t , U  \rangle =  \langle U , U \rangle = \frac{1}{\vert G \vert ^2 }. 
\end{equation}   
If $ \textbf{1}  : G \to U (1) $ denotes the trivial representation, then  (\ref{eq:fourier-properties}) and (\ref{eq:w-properties}) imply
\begin{equation}\label{eq:fourier-trivial}
    \widehat{W^{t}_{\alpha}} ( \textbf{1} ) = \hat{W}_{\alpha}^t(\textbf{1}) = \id _{\mathbb{C}} 
\end{equation}                               
for all  $t \in \mathbb{N}$. Then one has
\begin{equation}\label{eq:w-estimate}
    \begin{split}
\langle W^t_{\alpha} - U ,& W^t _{\alpha} - U \rangle  = \langle W^t_{\alpha}, W^t_{\alpha}\rangle  - \frac{1}{\vert G \vert ^2} \\
& = \frac{1}{\vert G \vert ^2} \sum _{\rho} d_{\rho} \, \tr ( \hat{W}_{\alpha}^{2t}(\rho) )  - \frac{1}{\vert G \vert ^2} \\
& = \frac{1}{\vert G \vert ^2} \sum _{\rho \neq 1} d_{\rho} \,  \tr ( \hat{W}_{\alpha}^{2t}(\rho) ) ,
    \end{split}
\end{equation}
where the last sum is taken over isomorphism classes of non-trivial irreducible unitary representations; the first equality uses (\ref{eq:w-u-product}), the second one follows from the Plancherel formula and (\ref{eq:w-properties}), and the third one uses (\ref{eq:fourier-trivial}). Finally, by the Cauchy-Schwarz inequality, (\ref{eq:w-estimate}) implies
\begin{equation}\label{eq:epsilon-estimate}
    \vert \varepsilon _{\alpha} (t) \vert ^2 \leq  \sum _{\rho \neq 1} d_{\rho}  \,\tr ( \hat{W}_{\alpha}^{2t}(\rho) ) .
\end{equation}
    When $\alpha  = 2 \vert S \vert$, we denote $W_{\alpha}^t$ by $W^t$, $W_{\alpha}$ by $W$, $\varepsilon_{\alpha}(t)$ by $\varepsilon (t)$, and $\tau_{\Gamma}^{\alpha}$ by $\tau _{\Gamma}$.

\section{Diameter bounds}\label{sec:diameter}

In this section we prove Theorems \ref{thm:diameter-effective}, \ref{thm:diameter-as}, and \ref{thm:diameter-effective-cayley}. For a group $G$ and a symmetric generating set $S\subset G$, we set $S_k : = S \cup S^2 \cup \ldots \cup S^{\left\lfloor  \frac{k+2}{3} \right\rfloor} $. It is straightforward to check that 
\begin{equation}\label{eq:sr-diam}
 \diam (\Gamma (G,S)) \leq \diam (\Gamma (G,S_k))  \left\lfloor \frac{k+2}{3} \right\rfloor .      
\end{equation}
The ensuing result \cite{behr} reduces the proof of Theorem \ref{thm:diameter-effective-cayley} to the case $k = 3$.
\begin{lemma}[Behr]\label{lem:behr}
\rm  Let $G$ a group, and $S\subset G$ a finite symmetric set of generators. Then for each $k \geq 3$, there is a surjective map 
\[\pi _1 (P_k (\Gamma (G,S))) \to \pi_1(P_3(\Gamma (G,S_k))).\] 
\end{lemma}

\begin{proof}
First notice that if an edge $\omega$ of $\Gamma (G, S_k) $ corresponds to an element of $S^m$ with $m \leq \lfloor \frac{k +2 }{3 }\rfloor$ then there is an endpoint-preserving homotopy in $P_3(\Gamma (G, S_k))$ taking $\omega$ to a concatenation of $m$ edges in $\Gamma (G, S)$ (see Figure \ref{fig:homot}). Since the fundamental group of a CW-complex is generated by the loops in its $1$-skeleton, the above observation implies that the inclusion
    \begin{equation}\label{eq:behr}
        \Gamma (G, S) \to P_3( \Gamma (G, S_k ) ) 
    \end{equation}                    
induces a surjective map at the level of fundamental groups. It remains to check that  (\ref{eq:behr}) extends to a continuous map $P_k(\Gamma (G,S)) \to P_3 (\Gamma (G,S_k) ) $.  This boils down to the fact that any word of length $\leq k$ representing the identity in $G$ using the elements of $S$ as letters, can be written as a concatenation of words of length $\leq 3$ representing the identity in $G$ using the elements of $S_k$ as letters (see \cite[Lemma 7.A.8]{cornulier-delaharpe} for further details).  
\begin{figure}[h!]
\centering
\psfrag{a}{$\omega$}
\psfrag{b}{$s_1$}
\psfrag{c}{$s_2$}
\psfrag{d}{$s_3$}
\psfrag{e}{$s_4$}
\psfrag{f}{$s_5$}
\psfrag{g}{$s_6$}
\includegraphics[scale = 0.5]{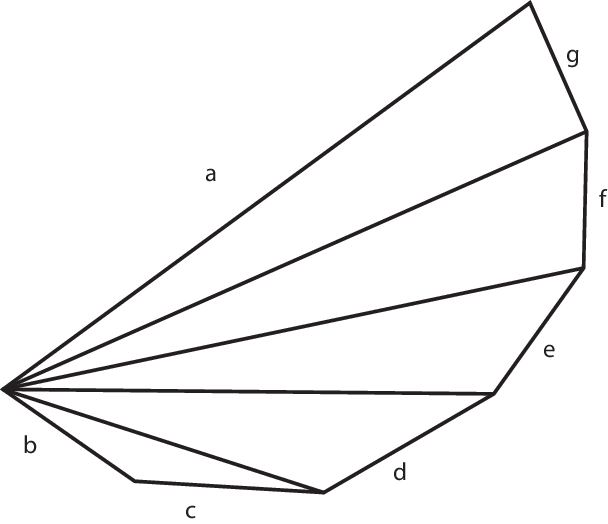}
\caption{If $\omega = s_1 \cdots s_m$, then the decomposition $w = (s_1 \cdots s_{m-1})(s_m)$ yields an endpoint preserving homotopy in $P_3(\Gamma (G,S_k))$ from $w$ to the concatenation of an edge corresponding to an element in $S^{m-1}$ and one in $S$. Proceeding inductively yields the desired homotopy. In the picture $m = 6$.}\label{fig:homot}
\end{figure}
\end{proof}

The following elementary observation will be required at the end of the proof of Theorem \ref{thm:diameter-effective-cayley}.

\begin{lemma}\label{lem:mv}
\rm Let $\Gamma$ be a connected graph,  $T \subset \Gamma$ a connected subgraph, and $C_1, \ldots , C_{\ell} \subset \Gamma $ the connected components of $\Gamma \backslash T$. Then for each $j \in \{ 1, \ldots , \ell \}$, the graph $\Gamma \backslash C_j$ is connected.
\end{lemma}

\begin{proof}
If the result is false, there are $w_0, w_1  \in \Gamma \backslash C_j$ such that any path connecting them passes through $C_j$. Since $\Gamma$ is connected, there is a path $[w_0 = v_0, v_1, \ldots , v_k = w_1]$, which by assumption passes through $C_j$. Let $i_1$ be the first index such that $v_{i_1 + 1} $ is in $C_j$ and let $i_2$ be the last index such that $v_{i_2 - 1}$ lies in $C_j$. Then both $v_{i_1}$ and $v_{i_2}$ lie in $T$, and by connectedness of $T$ there is a path $[v_{i_1} = a_0 , a_1, \ldots , a_{m} = v_{i_2} ]$ in $T$. Then it is easy to check that the path 
\[  [ w_0 = v_0, \ldots , v_{i_1} , a_1, \ldots , a_{m-1}, v_{i_2} , \ldots , v_k =w_1  ]           \]
does not intersect $C_j$, contradicting our assumption (see Figure \ref{fig:cj-con}). 
\begin{figure}[h!]
\centering
\psfrag{a}{$w_0 = v_0 $}
\psfrag{b}{$v_1$}
\psfrag{c}{$v_2$}
\psfrag{d}{$v_3$}
\psfrag{e}{$v_4$}
\psfrag{f}{$v_5$}
\psfrag{g}{$v_6$}
\psfrag{h}{$v_7$}
\psfrag{i}{$v_8$}
\psfrag{j}{$v_9$}
\psfrag{k}{$v_{10}$}
\psfrag{l}{$v_{11}$}
\psfrag{m}{$v_{12} = w_1$}
\psfrag{n}{$a_1$}
\psfrag{o}{$a_2$}
\psfrag{p}{$a_3$}
\psfrag{q}{$C_j$}
\includegraphics[scale = 1]{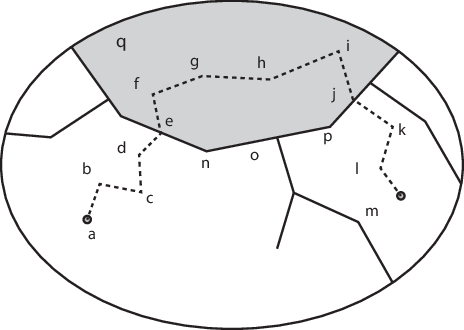}
\caption{The portion $[v_{i_1}, v_{i_1+1}, \ldots , v_{i_2-1}, v_{i_2}]$ can be replaced by $[a_0, \ldots, a_m]$. In the picture the shadowed region represents $C_j$, the solid line represents $T$, $i_1 = 4$, and $i_2=9$.}\label{fig:cj-con}
\end{figure}
\end{proof}

\begin{proof}[Proof of Theorem \ref{thm:diameter-effective-cayley}:]
By Lemma \ref{lem:behr} we have $\bet ( P_3(\Gamma (G, S_k)) )  = 0 $.  Combining this with  (\ref{eq:sr-diam}) we can assume $k = 3$ without loss of generality. 

Let $e \in G$ be the neutral element and $\Gamma  = \Gamma (G,S)$. Take $h\in G$ with $d_{\Gamma}(h,e)= m = \diam (\Gamma )$ and a minimizing path  $ [ e=g_0,g_1, \ldots , g_m=h ] $. For each $i$, set $\Sigma _i \subset \Gamma$ as the subgraph induced by the set of vertices $\{   g \in G \mid d_{\Gamma} (g,e) = i \}$ and let $T_i$ be the connected component of $\Sigma _i$ containing $g_i$. 

\begin{center}
Claim 1: For each $i_0 \in \{ 1, \ldots , m-1 \}$, the vertices $e$ and $h$ lie in 

distinct  connected components of $\Gamma \backslash T_{i_0}$. 
\end{center}
Let $Y_0 = Y_1  \subset  \Gamma $ be the subgraph induced by $\bigcup_{j=0}^{i_0-1} \Sigma _ j $, $Y_{1/4} \subset \Gamma$ the subgraph induced by $T_{i_0}$,   $Y_{1/2} \subset \Gamma$ the subgraph induced by $\bigcup _{j=i_0 + 1}^m \Sigma_j$, and $Y_{3/4} \subset \Gamma $ the subgraph induced by $ \Sigma_{i_0} \backslash T_{i_0} $. Since $d_{\Gamma} (\cdot , e )$ is $1$-Lipschitz, there are no edges between $Y_0$ and $Y_{1/2}$, and by the definition of $T_{i_0}$ there are no edges between $Y_{1/4}$ and $Y_{3/4}$. 

Then construct a map $\psi: \Gamma  \rightarrow  \mathbb{R} / \mathbb{Z} $  that restricted to $Y_s$ equals $s$ for $s \in \{ 0, 1/4, 1/2, 3/4 \}$, and takes all edges joining $Y_0$ with $Y_{1/4}$ to the interval $[0, 1/4]$, doing the same for the intervals $[1/4, 1/2]$, $[1/2, 3/4]$, and $[3/4, 1]$. 

By construction, $\psi$ sends each edge of $\Gamma $ to either a point or an interval of length $1/4$ in $ \mathbb{R} / \mathbb{Z} $. Also recall that each $2$-cell $\alpha$ of $P_3(\Gamma )$ is attached to $\Gamma$ via a simple loop $\partial \alpha$ of length $\leq 3$. Therefore, $\psi (\partial \alpha ) \subset  \mathbb{R} / \mathbb{Z} $ is a loop of length $\leq 3/4$ hence nullhomotopic for each $\alpha$, and $\psi$ extends to a map $\Psi : P_3(\Gamma ) \to  \mathbb{R} / \mathbb{Z} $ (see Figure \ref{fig:y1234}).

\begin{figure}[h!]
\centering
\psfrag{a}{$Y_0 = Y_1$}
\psfrag{b}{$Y_{1/4}$}
\psfrag{c}{$Y_{1/2}$}
\psfrag{d}{$Y_{3/4}$}
\psfrag{e}{$0=1$}
\psfrag{f}{$1/4$}
\psfrag{g}{$1/2$}
\psfrag{h}{$3/4$}
\psfrag{i}{$\Psi$}
\psfrag{j}{$\Phi$}
\includegraphics[scale = 1.2]{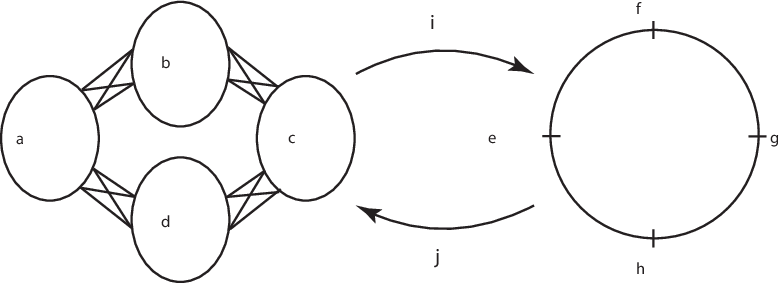}
\caption{Each edge of $\Gamma$ is sent via $\psi$ to either a point or an interval of length $\leq 1/4$, so $\psi$ extends to a map $\Psi : P_3(\Gamma) \to \mathbb{R}/\mathbb{Z}$.}\label{fig:y1234}
\end{figure}

 Assume the claim is false and take a minimizing path $[ e=g^{\prime}_0,g^{\prime}_1, \ldots , g^{\prime}_{m^{\prime}}=h ] $  in $\Gamma \backslash T_{i_0} $. Consider the map $\Phi :  \mathbb{R}/\mathbb{Z} \rightarrow P_3 (\Gamma  )$ that sends the interval $[0, 1/2] $ to the path $[g_0,g_1, \ldots ,g_m]$ and $[1/2, 1]$ to  the path $[g^{\prime}_{m^{\prime}} , \ldots , g^{\prime}_1 , g^{\prime}_0]$. 

Since $\Psi \circ \Phi (0) = 0$, $\Psi \circ \Phi (1/2) = 1/2$,  $\Psi \circ \Phi \vert _{[0, 1/2]}$ misses $3/4$, and $\Psi \circ \Phi \vert _{[1/2, 1]} $ misses $1/4$, the composition $\Psi \circ \Phi$ is homotopic to the identity in $ \mathbb{R} / \mathbb{Z} $, meaning that the induced map $\Psi_{\ast} : \pi _1(P_3(\Gamma)) \to \pi _1( \mathbb{R} / \mathbb{Z} )= \mathbb{Z}$ is surjective, contradicting the hypothesis $\bet (P_3(\Gamma)) =0$. This finishes the proof of the claim.

\begin{center}
Claim 2: For each $i_0 \in \{ 1, \ldots, m-1 \}$, either 

$g_{i_0}^{-1}T_{i_0}$ or $hg_{i_0}^{-1}T_{i_0}$ intersect $T_{i_0}$.
\end{center}
Let $C_1, \ldots , C_{\ell} \subset \Gamma $ denote the connected components of $\Gamma \backslash T_{i_0}$, with $e \in C_1$, $h \in C_2$. Observe that $g_{i_0}^{-1}C_1, \ldots , g_{i_0}^{-1}C_{\ell}$ are the connected components of $\Gamma \backslash g_{i_0}^{-1}T_{i_0}$, and  $hg_{i_0}^{-1}C_1, \ldots , $ $hg_{i_0}^{-1}C_{\ell}$ are the connected components of $\Gamma \backslash hg_{i_0}^{-1}T_{i_0}$. 

Assume $T_{i_0}  \cap  g_{i_0}^{-1}T_{i_0} = \emptyset $.  Since  $g_{i_0}^{-1}T_{i_0}$ is connected and contains $e$, it is contained in $ C_1 $. By Lemma \ref{lem:mv}, $\Gamma \backslash C_1 = T _{i_0}\cup C_2 \cup \ldots  \cup C_{\ell}$ is connected, and since it doesn't intersect $g_{i_0}^{-1}T_{i_0}$, it is contained in $g_{i_0}^{-1}C_{j_1}$  for some $j_1$. Therefore 
\[   \vert C_j \vert < \vert C_{j_1} \vert \text{ for all } j \neq 1 .                \]
This is only possible if $j_1 = 1$, and in particular we have 
\begin{equation}\label{eq:c1-c2}
    \vert C_2 \vert < \vert C_1 \vert  .    
\end{equation}      
Similarly, if $T_{i_0}  \cap  hg_{i_0}^{-1}T_{i_0} = \emptyset $, then $hg_{i_0}^{-1}T_{i_0}$ is contained in $C_2$. By Lemma \ref{lem:mv}, $\Gamma \backslash C_2 = T_{i_0} \cup C_1 \cup C_3 \cup \ldots  \cup C_{\ell}$ is connected, and since it doesn't intersect $hg_{i_0}^{-1}T_{i_0}$, it is  contained in $h g_{i_0}^{-1}C_{j_2}$ for some $j_2$, meaning that 
\[         \vert C_j \vert < \vert C_{j_2} \vert \text{ for all } j \neq 2 .                                            \]
This implies that $j_2 = 2$ and 
\begin{equation}\label{eq:c2-c1}
  \vert C_1 \vert < \vert C_2 \vert  .
\end{equation}
Assuming the claim is false, both (\ref{eq:c1-c2}) and (\ref{eq:c2-c1}) would hold; a contradiction.

From the second claim,  we deduce  $\textnormal{diam} (T_i) \geq \min \{   i, m-i \}$ for each $i \in \{ 1 , \ldots , m -1 \}$, and  since the $T_i$'s are disjoint, we conclude that 
 \[ \vert G \vert \geq \Sigma_{i=0}^m \vert T_i \vert \geq  \Sigma_{i=0}^m \min   \{ i+1, m-i+1  \}.  \]
  This implies that $m \leq \sqrt{4\vert G \vert +1} -2$, which is the required inequality. \end{proof}

\begin{proof}[Proof of Theorem \ref{thm:diameter-effective}:]
 Let $p \in X$ and define $S$ as in Proposition \ref{prop:sml} with $\delta = 0 $. By Proposition \ref{prop:cover}, the first Betti number of $P_3(\Gamma (G, S))$ vanishes. For $x, y \in X$,  take $g_1, g_2 \in G $ with $d_X(g_1 p, x), d_X(g_2 p, y ) \leq \diam (X/G)$. Then 
 \begin{align*}
     d_X(x,y) &\leq 2 \cdot \diam (X/G) + d_X(g_1 p , g_2 p) \\
     &\leq  2 \cdot \diam (X/ G) [1 + d_{\Gamma}(g_1, g_2)]\\
     & \leq 2 \cdot \diam (X/G) [ \sqrt{4 \vert G \vert + 1} - 1 ],
  \end{align*}
 where the second inequality follows from \v{S}varc--Milnor Lemma and the third one from Theorem \ref{thm:diameter-effective-cayley}. The result follows since $\sqrt{4 \vert G \vert + 1} \leq 2 \sqrt{\vert G \vert } + 1$. 
\end{proof}

\begin{proof}[Proof of Theorem \ref{thm:diameter-as}:]
Pick $p_n \in X_n$, set $S_n : = \{ g \in G_n \vert d(gp_n, p_n) \leq 2 \cdot \diam (X_n / G_n ) \} $, and let $\Gamma _n : = \Gamma (G_n, S_n)$.  If the result fails, there is $\varepsilon > 0$ such that after taking a subsequence one has
\[       \vert G_n \vert ^{\varepsilon} = O \left( \frac{\diam (X_n)}{\diam (X_n / G_n)} \right)    .      \]
By the \v{S}varc--Milnor Lemma, this would imply  $ \vert G_n \vert ^{\varepsilon} = O \left( \diam (\Gamma_n ) \right) .  $  Then by \cite[Theorem 1]{benjamini-finucane-tessera}, after further taking a subsequence, $\Gamma_n / \diam (\Gamma_n)$ converges to an $m$-dimensional torus $X$. By Remark \ref{rem:p-metric}, the sequence $P_3(\Gamma_n) / \diam (P_3(\Gamma_n))$ also converges to $X$ and by \cite[Theorem 2.1]{sormani-wei}  there are surjective morphisms $\pi_1 (P_3(\Gamma_n)) \to \pi_1 (X) =  \mathbb{Z}^m$ for large enough $n$, contradicting Proposition \ref{prop:cover}. 
\end{proof}

\section{Fourier analysis in abelian groups}\label{sec:abelian}

In this section we prove Theorem \ref{thm:k-spectrum} and Corollaries \ref{cor:mix} and \ref{cor:ricci}.  Let $k$, $G$, $S$ be as in the statement of Theorem \ref{thm:k-spectrum} and let $\Gamma : = \Gamma (G,S) $.  Assuming the estimate (\ref{eq:K}) fails to hold, there is an irreducible non-trivial unitary representation $\rho : G \to U(m)$ and $x \in \mathbb{S}^{2m-1}$ with 
\begin{equation}\label{eq:no-K}
    \sup_{s \in S} d(\rho (s) x , x) < 2 \cdot \sin ( \pi / k )  .     
\end{equation}      
By Proposition \ref{prop:abelian-1} we have $m = 1$, and since the metric on $\mathbb{S}^1 $ is bi-invariant, we can assume $ x = 1$. Let $\psi : \Gamma \to \mathbb{S}^{1}$ be the map that restricted to $G$ coincides with $\rho$, and restricted to an edge  $[g,h] \subset \Gamma $ is a minimizing geodesic from $\rho (g)$ to $\rho (h)$. 

If $g,h \in G$ are such that $g = hs$ for some $s \in S$, then (\ref{eq:no-K}) implies that the angle between $\rho (g) $ and $\rho (h)$ is less than  $2 \pi / k $. Hence, for any simple loop of length $\leq k$ in $\Gamma$, its image under $\psi$ has length less than $2 \pi$ and is contractible. Therefore, $\psi$ extends to a map $\Psi : P_k (\Gamma) \to \mathbb{S}^1$. 

Since $\rho$ is non-trivial, there is $s \in S$ with $\rho (s) \neq 1$. Then the image under $\Psi$ of the loop $[ e, s, s^2, \ldots , s^{\vert G \vert } =e]$ is a loop in $\mathbb{S}^1$ that winds around at least once; counterclockwise if $\real (\rho (s))>0$ and clockwise if $\real (\rho (s) ) < 0$. This means the map $ \Psi _{\ast }:\pi_1 (P_k (\Gamma ) ) \to \pi_1(\mathbb{S}^1) = \mathbb{Z}$ is non-trivial, contradicting the assumption $\bet (P_k(\Gamma)) = 0$. This finishes the proof of (\ref{eq:K}).

We now proceed to prove  (\ref{eq:lambda}). Notice that if we simply apply  (\ref{eq:isoperimetrics}) naively to (\ref{eq:K}) we would get a weaker result. Let $ - \lambda \in \sigma (\Delta )\backslash \{ 0 \} $, and  $E_{\lambda} \leq \mathbb{C} [G]$  the corresponding eigenspace. By Lemma \ref{lem:x-irreducible}, there is $x \in E_{\lambda} $  with $\vert x \vert = 1$ such that the span of its orbit is an  irreducible representation $\rho : G \to \mathbb{S}^1$. Since $\lambda \neq 0$, it follows that $\rho $ is non-trivial. 
\begin{center}
    Case 1: $\rho (s) \neq -1$ for all $s \in S$.
\end{center}
By (\ref{eq:K}), there is $t \in S$ with $d(\rho(t), 1) \geq 2 \cdot \sin (\pi / k)$. This implies $\real (\rho (t) )  = \real (\rho (t^{-1})) \leq \cos ( 2 \pi / k )$. Since $\rho (t) \neq -1$, one has $t \neq t^{-1}$ so
\begin{align*}
     \lambda =  \langle  - \Delta \,  x, x \rangle &  = \sum_{s \in S} [ 1- \rho (s)    ]\\
    &  =  \vert S \vert  -  \sum_{s \in S} \real (\rho (s))\\
    & \geq 2 - 2 \cos (2 \pi / k ),
 \end{align*}
 where in the last line we used that $\real (\rho (s)) < 1$ for all $s \in S \backslash \{ t, t^{-1} \} $.
\begin{center}
    Case 2: $k = 3 $ and $\rho (t) = -1$ for some $t \in S$. 
\end{center}
We claim there are $s_1, s_2 \in S \backslash \{ e , t \}$ such that $s_1s_2 = t$. Assume otherwise; then the edge $[e,t ] \subset \Gamma $ does not belong to a $2$-cell of $P_3(\Gamma)$. Let $x$ be the midpoint of the edge $[e,t]$, $A := [e,t]$, and $B := P_3(\Gamma) \backslash \{ x \}$. Then $A \cap B = [e , x) \cup (x , t ]$ and by Lemma \ref{lem:t-not-separating}, $B$ is connected,  so the portion of the Mayer--Vietoris sequence (with real coefficients) 
\[    H_1(P_3(\Gamma)) \to H_0(A \cap B ) \to H_0(A) \oplus H_0 (B) \to H_0 (P_3(\Gamma)) \to 0   \]
yields the exact sequence (using $\bet (P_3(\Gamma)) = 0$)
\[    0 \to \mathbb{R} \oplus \mathbb{R} \to \mathbb{R} \oplus \mathbb{R} \to \mathbb{R} \to 0 .                          \]
This is impossible by dimension counting. Then there are $s_1, s_2 \in S \backslash \{ e, t \}$ (not necessarily distinct) with $ s_1 s_2= t$.  Since $\rho (s_1)\rho(s_2) = -1$, without loss of generality we can assume $\real (\rho (s_1)) \leq 0 $. Then 
\begin{align*}
    \lambda & =   \vert S \vert  -  \sum_{s \in S} \real (\rho (s)) \geq 2 - \real (\rho (t) + \rho (s_1))   \\
    &\geq 3 = 2 - 2 \cos (2 \pi /3 ).
\end{align*}      
\begin{center}
    Case 3: $k \geq 4$ and $\rho (t) = -1$ for some $t \in S$. 
\end{center}
Directly compute
\[     \lambda  = \vert S \vert  -  \sum_{s \in S} \real (\rho (s)) \geq 2 \geq 2 - 2 \cos (2 \pi / k ). \]
This, together with (\ref{eq:isoperimetrics}), finishes the proof of (\ref{eq:lambda}). 

For notational convenience, we set 
\[ \xi  _k : = 1 - \cos (2 \pi / k )\, \text{  and }\, \alpha_k : =  \vert S \vert + \xi _k. \]
To prove (\ref{eq:diameter}), we look at the random walk $W_{\alpha_k}^t$ in $\Gamma$.
\begin{lemma}\label{lem:w-hat}
    \rm  For any non-trivial irreducible unitary representation $\rho : G \to \mathbb{S}^1$,
    \[ \hat{W}_{\alpha_{k}} (\rho) \in \left[ -1 +   \frac{2 \xi_k}{\alpha_k}   , 1 - \frac{2 \xi_k}{\alpha_k}    \right] .\]
\end{lemma}
\begin{proof}
This is a direct computation using (\ref{eq:w1-expression}). On one hand we have
\[
    \hat{W}_{\alpha_{k}} (\rho )    =  1 + \frac{1}{\alpha_k} \sum _{s \in S} [ \rho (s) - 1 ]  \leq   1 - \frac{2 \xi_k}{\alpha_k}, 
\] 
where we first used the identity  $\hat{\delta}_s (\rho) = \rho (s)$ and then the estimate $\lambda_1 \geq 2 \xi_k$.  For the other inequality, notice that $\alpha_k - 2 \vert S \vert  = 2 \xi _k - \alpha_k$, then
\[    \hat{W}_{\alpha_{k}} (\rho )    =   1 + \frac{1}{\alpha_k} \sum _{s \in S} [ \rho (s) - 1 ]   \geq \frac{\alpha_k - 2 \vert S \vert  }{\alpha_k  }     = \frac{ 2\xi_k  }{\alpha_k}  - 1  ,
\] 
where we used $\real (\rho (s)) \geq -1 $ for all $s \in S$ in the inequality.  
\end{proof}
Then by (\ref{eq:w-estimate}),  for $t \in \mathbb{N}$ we have
\[
   \langle W^t_{\alpha_k} - U , W^t_{\alpha_k} - U \rangle  =  \frac{1}{\vert G \vert ^2 } \sum_{\rho \neq 1} d_{\rho} \, \tr \left( \hat{W}^{2t}_{\alpha_k} \right)  <  \frac{1}{\vert G \vert  }\left\vert 1 - \frac{2 \xi_k}{\alpha_k} \right\vert ^{2t} ,
\]
where we used (\ref{eq:sum-dimensions}) for the inequality. Now assume $t \leq \diam (\Gamma)$. Since $W^t_{\alpha_k}$ is supported in the ball of radius $t$ around $e$, then the left hand side of the equation is at least $ \frac{1}{\vert G \vert ^3}$. Hence by taking logarithm and using the identity $\log (1 + u ) \leq u$ we get
\[  -  2 \log \vert G \vert  <  2t \log     \left\vert 1 - \frac{2 \xi_k}{\alpha_k} \right\vert  \leq  -  \frac{4t \xi _k}{\alpha_k} .         \]
Rearranging terms, 
\[    t < \frac{  \alpha_k \log \vert G \vert  }{ 2 \xi_k }.    \]
This implies 
\[           \diam (\Gamma) \leq  \frac{\alpha_k  }{2\xi_k } \log  \vert G \vert  + 1.                                 \]
This concludes the proof of Theorem \ref{thm:k-spectrum}. Notice that if  $k = 3$, then (\ref{eq:diameter}) simplifies to
\begin{equation}\label{eq:diameter-3}
     \diam (\Gamma) \leq  \left[ \frac{\vert S \vert  }{3 } + \frac{1}{2} \right] \log  \vert G \vert  + 1.         
\end{equation}                      
In order to prove Corollary \ref{cor:mix}, we need to establish an analogue of Lemma \ref{lem:w-hat}.
\begin{lemma}\label{lem:w-walk}
    \rm  Under the hypothesis of Corollary \ref{cor:mix}, for any non-trivial irreducible unitary representation $\rho : G \to \mathbb{S}^1$ one has
    \[ \hat{W} (\rho) \in \left[  0  , 1 - \frac{\xi_k}{\vert S \vert }    \right] .\]
\end{lemma}
\begin{proof}
This is again a direct computation. Using $\lambda _ 1 \geq 2 \xi_k$ we get
\[
    \hat{W} (\rho )    =  1 + \frac{1}{2 \vert S \vert } \sum _{s \in S} [ \rho (s) - 1 ]  \leq   1 - \frac{\xi_k}{\vert S \vert }. 
\] 
On the other hand, simply using $\real (\rho (s)) \geq -1$ for all $s \in S$ we conclude 
\[    \hat{W} (\rho )    =   1 + \frac{1}{2 \vert S \vert } \sum _{s \in S} [ \rho (s) - 1 ]   \geq 0  .
\] 
\end{proof}
\begin{proof}[Proof of  Corollary \ref{cor:mix}:]
By (\ref{eq:epsilon-estimate}), using  Lemma \ref{lem:w-walk} and (\ref{eq:sum-dimensions}), we have
\[ 
       \varepsilon (t) ^2  \leq \sum_{\rho \neq 1} d_{\rho} \, \tr (\hat{W}^{2t})  < \vert G \vert \left\vert 1 - \frac{\xi_k}{\vert S \vert } \right\vert ^{2t}  
\] 
for $t \in \mathbb{N}$. If $\varepsilon (t) \geq c \in [0, 2] $, then taking logarithms as above we get
\[          2 \log (c)  < \log \vert G \vert  -   \frac{2 t \xi_k}{\vert S \vert } .                                  \]
Rearranging terms we get
\[       t < \frac{ \vert S \vert }{ 2   \xi_k }[ \log \vert G \vert - 2 \log (c) ].           \]
Since $\xi_k \geq 16/k^2$ the result follows. 
\end{proof}

\begin{proof}[Proof of Corollary \ref{cor:ricci}:]
Take  a point $p \in M $ with injectivity radius $\geq 2r_0$,  $\tilde{p} \in \tilde{M}$ in its preimage, and set 
\[     S : = \{  g \in \pi_1(M) \backslash \{ e \} \, \vert \, d(g \tilde{p}, \tilde{p}) \leq 2 D  \}  .  \]
By the injectivity radius condition, for $g,h \in \pi_1(M)$ distinct, the balls $B(g\tilde{p}, r_0)$ and $B(h\tilde{p}, r_0)$ are isometric and disjoint. Since $g \in S \cup \{ e \}$ implies $B(g\tilde{p}, r_0) \subset B(\tilde{p}, 2D + r_0)$, we have  
\[    \vert S \vert +1 \leq \frac{ \vol ( B( \tilde{p},  2 D + r_0 ) ) }{ \vol (B (\tilde{p}, r_0)  )}                      .      \]
By the Bishop--Gromov inequality \cite[Section 11.10]{bishop-crittenden}, the right hand side of the equation is less or equal than $v _{n}^{\kappa} (2D + r_0) / v_n^{\kappa}(r_0)$. By Proposition \ref{prop:cover}, $P_3(\Gamma (\pi_1(M), S))$ is simply connected, so (\ref{eq:diameter-3}) holds;
\begin{equation}\label{eq:ricci-1}
     \diam (\Gamma (\pi_1(M), S)) \leq \left[    \frac{v _{n}^{\kappa} (2D + r_0)}{3v_n^{\kappa}(r_0)}  + \frac{1}{6}          \right] \log \vert \pi_1(M) \vert  + 1 .
\end{equation}                        
Arguing as in the proof of Theorem \ref{thm:diameter-effective}, we have
\begin{equation}\label{eq:ricci-2}
  \frac{\diam (\tilde{M})}{\diam(M)} \leq 2 + 2 \cdot  \diam (\Gamma (\pi_1 (M), S))  .                      
\end{equation}
Combining (\ref{eq:ricci-1}) and (\ref{eq:ricci-2}) the result follows.
\end{proof}

\section*{Acknowledgements}
The author would like to express his thanks to Cameron Rudd and Andrew Ng for suggesting Example \ref{exam:hyperbolic}, and  Jacob Bradd, Ana C. Ch\'avez-C\'aliz, and Anton Petrunin for lengthy discussions and helpful comments about previous versions of this paper. 

He is also grateful to Sebastian Gouzel for pointing out the ideas in \cite{brooks-i}, \cite{brooks-ii}, \cite{burger}, \cite{magee} to him, stimulating his interest in Theorem \ref{thm:k-spectrum} even further, and an anonymous reviewer whose contributions have significantly improved not only the presentation of this paper but the quality of the results.

Starting from August 2022, the author holds a Postdoctoral Fellowship at the Max Planck Institute for Mathematics at Bonn. 


\begin{thebibliography}{}
%
%


\bibitem{ash}  Ash, Robert B. Probability and measure theory. Second edition. With contributions by Catherine Dol\'eans-Dade Harcourt/Academic Press, Burlington, MA, 2000. xii+516 pp. \textbf{ISBN}:0-12-065202-1




\bibitem{babai-seress} Babai, L\'aszl\'o; Seress, \'Akos. \textit{On the diameter of permutation groups}. European J. Combin. \textbf{13} (1992), no.4, 231-243.

\bibitem{behr} Behr, Helmut.  \textit{Endliche Erzeugbarkeit arithmetischer Gruppen \"uber Funktionenk\"orpern}. Invent. Math. \textbf{7} (1969), 1-32.

\bibitem{bekka-delaharpe-valette} Bekka, Bachir; de la Harpe, Pierre; Valette, Alain. Kazhdan's property (T). New Math. Monogr., 11. Cambridge University Press, Cambridge, 2008. xiv+472 pp. \textbf{ISBN}:978-0-521-88720-5

\bibitem{benjamini-finucane-tessera}  Benjamini, Itai; Finucane, Hilary; Tessera, Romain
\textit{On the scaling limit of finite vertex transitive graphs with large diameter}. Combinatorica \textbf{37} (2017), no.3, 333-374.

\bibitem{bishop-crittenden} Bishop, Richard L.; Crittenden, Richard J. Geometry of manifolds. Reprint of the 1964 original. 
AMS Chelsea Publishing, Providence, RI, 2001. xii+273 pp. \textbf{ISBN}:0-8218-2923-8


\bibitem{boston-ellenberg} Boston, Nigel; Ellenberg, Jordan S. \textit{Pro-p groups and towers of rational homology spheres}. Geom. Topol. \textbf{10} (2006), 331-334.

\bibitem{breuillard-green-tao} Breuillard, Emmanuel; Green, Ben; Tao, Terence. \textit{The structure of approximate groups}. Publ. Math. Inst. Hautes \'Etudes Sci. \textbf{116} (2012), 115-221.


\bibitem{breuillard-tointon} Breuillard, Emmanuel; Tointon, Matthew C. H. \textit{Nilprogressions and groups with moderate growth}. Adv. Math. \textbf{289} (2016), 1008-1055.

\bibitem{brooks-i}  Brooks, Robert. \textit{The fundamental group and the spectrum of the Laplacian}. Comment. Math. Helv. \textbf{56} (1981), no.4, 581-598.

\bibitem{brooks-ii}  Brooks, Robert. \textit{The spectral geometry of a tower of coverings}. J. Differential Geom. \textbf{23} (1986), no.1, 97-107.

\bibitem{brooks-iii} Brooks, Robert. \textit{Some remarks on volume and diameter of Riemannian manifolds}. J. Differential Geom. \textbf{27} (1988), no.1, 81-86.

\bibitem{burger} Burger, Marc. \textit{Spectre du laplacien, graphes et topologie de Fell}. Comment. Math. Helv. \textbf{63} (1988), no.2, 226-252

\bibitem{buser} Buser, Peter. \textit{A note on the isoperimetric constant}. Ann. Sci. \'Ecole Norm. Sup. (4) \textbf{15} (1982), no.2, 213-230.

\bibitem{calegari-dunfield} Calegari, Frank; Dunfield, Nathan M.
\textit{Automorphic forms and rational homology 3-spheres}. Geom. Topol. \textbf{10} (2006), 295-329.

\bibitem{cornulier-delaharpe} Cornulier, Yves; de la Harpe, Pierre. Metric geometry of locally compact groups. EMS Tracts Math., 25 European Mathematical Society (EMS), Z\"urich, 2016. viii+235 pp. \textbf{ISBN}:978-3-03719-166-8



\bibitem{delaharpe} de la Harpe, Pierre. Topics in geometric group theory. Chicago Lectures in Math. University of Chicago Press, Chicago, IL, 2000. vi+310 pp. \textbf{ISBN}:0-226-31719-6
 \textbf{ISBN}:0-226-31721-8

\bibitem{delasalle-tessera} de la Salle, Mikael; Tessera, Romain.
\textit{Characterizing a vertex-transitive graph by a large ball}. J. Topol. \textbf{12} (2019), no.3, 705-743.

\bibitem{diaconis} Diaconis, Persi. Group representations in probability and statistics. IMS Lecture Notes Monogr. Ser., 11. Institute of Mathematical Statistics, Hayward, CA, 1988. vi+198 pp. \textbf{ISBN}:0-940600-14-5

\bibitem{eberhard-jezernik} Eberhard, Sean; Jezernik, Urban. \textit{Babai's conjecture for high-rank classical groups with random generators}. Invent. Math. \textbf{227} (2022), no.1, 149-210.

\bibitem{gelander}  Gelander, Tsachik. \textit{Limits of finite homogeneous metric spaces}. Enseign. Math. (2) \textbf{59} (2013), no.1-2, 195-206.

\bibitem{gorodski-lange-lytchak-mendes}  Gorodski, Claudio; Lange, Christian; Lytchak, Alexander; Mendes, Ricardo A. E. \textit{A diameter gap for quotients of the unit sphere}. J. Eur. Math. Soc. (JEMS) \textbf{25} (2023), no.9, 3767-3793. 

\bibitem{gromov} Gromov, Misha. Metric structures for Riemannian and non-Riemannian spaces.
Based on the 1981 French original. With appendices by M. Katz, P. Pansu and S. Semmes. Translated from the French by Sean Michael Bates. Reprint of the 2001 English edition. Mod. Birkh\"auser Class. Birkh\"auser Boston, Inc., Boston, MA, 2007. xx+585 pp. \textbf{ISBN}:978-0-8176-4582-3. \textbf{ISBN}:0-8176-4582-9

\bibitem{halasi} Halasi, Zolt\'an. \textit{Diameter of Cayley graphs of SL(n,p) with generating sets containing a transvection}. J. Algebra \textbf{569} (2021), 195-219.

\bibitem{ledoux} Ledoux, Michel. The concentration of measure phenomenon. Math. Surveys Monogr., 89. American Mathematical Society, Providence, RI, 2001. x+181 pp. \textbf{ISBN}:0-8218-2864-9.

\bibitem{lubotzky} Lubotzky, Alexander. Discrete groups, expanding graphs and invariant measures. With an appendix by Jonathan D. Rogawski. Progr. Math., 125. Birkh\"auser Verlag, Basel, 1994. xii+195 pp. \textbf{ISBN}:3-7643-5075-X


\bibitem{lubotzky-i} Lubotzky, Alexander. \textit{Eigenvalues of the Laplacian, the first Betti number and the congruence subgroup problem}. Ann. of Math. (2) \textbf{144} (1996), no.2, 441-452.


\bibitem{magee} Magee, Michael. \textit{On Selberg's eigenvalue conjecture for moduli spaces of abelian differentials}. Compos. Math. \textbf{155} (2019), no.12, 2354-2398.


\bibitem{mantuano} Mantuano, Tatiana. \textit{Discretization of compact Riemannian manifolds applied to the spectrum of Laplacian}. Ann. Global Anal. Geom. \textbf{27} (2005), no.1, 33-46.


\bibitem{martelli} Martelli, B.: An introduction to geometric topology. arXiv preprint arXiv:1610.02592 (2016)

\bibitem{petrunin-general} A. Petrunin.: Diameter of $m$-fold cover. https://mathoverflow.net/questions/7732/ diameter-of-m-fold-cover


\bibitem{petrunin-universal} A. Petrunin.: Diameter of universal cover. https://mathoverflow.net/questions/8534/ diameter-of-universal-cover.

\bibitem{serre} Serre, Jean-Pierre. Linear representations of finite groups. Translated from the second French edition by Leonard L. Scott. Grad. Texts in Math., Vol. 42. Springer-Verlag, New York-Heidelberg, 1977. x+170 pp. \textbf{ISBN}:0-387-90190-6


\bibitem{sormani-wei} Sormani, Christina; Wei, Guofang. \textit{Hausdorff convergence and universal covers}. Trans. Amer. Math. Soc. \textbf{353} (2001), no.9, 3585-3602.

\bibitem{todhunter} I. Todhunter: Spherical Trigonometry
I. Todhunter. Hawk Press. \textbf{ISBN/EAN}: 9789393971593

\bibitem{thurston} Thurston, William P. Three-dimensional geometry and topology. Vol. 1. Edited by Silvio Levy
Princeton Math. Ser., 35. Princeton University Press, Princeton, NJ, 1997. x+311 pp. \textbf{ISBN}:0-691-08304-5

\bibitem{turing} Turing, A. M. \textit{Finite approximations to Lie groups}. Ann. of Math. (2) \textbf{39} (1938), no.1, 105-111.

\bibitem{zamora} Zamora, S. \textit{Limits of almost homogeneous spaces and their fundamental groups}. to appear in Groups, Geometry and Dynamics. 




\end{thebibliography}

%

\end{document}